\documentstyle[12pt]{article}
\begin{document}
\def\cross{{>\!\!\!\triangleleft}}
\begin{center}
{\bf IMPLICATIONS OF THE HOPF ALGEBRA PROPERTIES \\
OF NONCOMMUTATIVE DIFFERENTIAL CALCULI}
\end{center}
\begin{center}
{\sc A. A. Vladimirov}
\end{center}
\begin{center}
{\em Bogoliubov Laboratory of Theoretical Physics,
Joint Institute for Nuclear Research, \\
Dubna, Moscow region 141980, Russia}
\end{center}

\vspace{.5cm}

{\small We define a noncommutative algebra of four basic
objects within a differential calculus on quantum groups -- functions,
1-forms, Lie derivatives, and inner derivations -- as the cross-product
algebra associated with Woronowicz's (differential) algebra of 
functions and forms. This definition properly takes into account the
Hopf algebra structure of the Woronowicz calculus. It also provides
a direct proof of the Cartan identity.}

\vspace{1cm}

In the framework of Woronowicz's noncommutative differential calculus
~\cite{Wo,AC} one deals with a differential complex
\begin{equation}
A\stackrel{d}{\longrightarrow}\Gamma\stackrel{d}{\longrightarrow}
\Gamma^2\longrightarrow\ldots \ \ , \label{1}
\end{equation}
where $A$ is a Hopf algebra (of functions on a quantum group),
 $\Gamma$ is a bicovariant $A$-bimodule,
$\Gamma^2\equiv \Gamma\wedge\Gamma$ its second (wedge) power,
and so on. Brzezinski~\cite{Br} has shown that 
\begin{equation}
\Gamma^{\wedge}\doteq A\oplus\Gamma\oplus\Gamma^2\oplus\ldots \label{2}
\end{equation}
also becomes a (graded) Hopf algebra with respect to (wedge)
multiplication and natural definitions of coproduct and antipode.
In the present note, we want to demonstrate how this Hopf structure
can be used to construct a noncommutative algebra containing 
4 types of elements: functions (on a quantum group), differential
1-forms, Lie derivatives, and inner derivations. Similar algebras
were introduced and studied by several 
authors~\cite{SWZ1,SWZ2,Is-4}.  Probably, the most close to ours is 
the approach by P.Schupp~\cite{Sch-Vare}. However, some of our 
results and, especially, starting points appear different. So, we 
propose the principles and recipes described below as fully 
Hopf-algebra motivated (and, we believe, natural) new approach in this 
field.

Let us recall the basic definitions of Woronowicz's differential
calculus:
\begin{equation}
 \omega ^ia=(f_{j}^{i}\triangleright a)\omega ^j\,, \ \ \ 
\Delta(\omega ^i)=1\otimes \omega ^i + \omega ^j\otimes r_{j}^{i} \,,
 \ \ \ da=(\chi _i\triangleright a)\omega ^i \,. 
\label{3}
\end{equation}
Here $a$ is arbitrary element of $A$, \ $r_{j}^{i}\in A$, \ $\chi _i$
and $f_{j}^{i}$ belong to $A^*$, and left-invariant 1-forms $\omega ^i
\in \Gamma$. The consistency conditions of the
calculus are
\begin{equation}
 \Delta(r_{j}^{i})=r_{j}^{k}\otimes r_{k}^{i}\,, \ \
\Delta(f_{j}^{i})=f_{k}^{i}\otimes f_{j}^{k}\,, \ \ 
(f_{i}^{j}\triangleright a)r_{k}^{i}=r_{i}^{j}(a\triangleleft 
f_{k}^{i}) \,, \ \ \varepsilon(r_{j}^{i})=\varepsilon (f_{j}^{i})=
\delta _{j}^{i}\,, \label{4} 
\end{equation} 
\begin{equation} \Delta(\chi _i)=\chi 
_j\otimes f_{i}^{j}+1\otimes \chi _i\,, \ \ \ a\triangleleft\chi 
_i=(\chi _j\triangleright a)r_{i}^{j}\,, \ \ \ 
\varepsilon (\chi_i)=0\,. \label{4a} 
\end{equation} 
We use the notation 
\begin{equation}
x\triangleright a\doteq a_{(1)}<\!x,a_{(2)}\!>\,, \ \ \  
a\triangleleft x\doteq a_{(2)}<\!x,a_{(1)}\!>\,, \label{4.5}
\end{equation}
where $\Delta(a)\equiv a_{(1)}\otimes a_{(2)}$, \ for left and right
actions, respectively, of some $x\in A^*$ upon $A$. Both these actions
are `covariant'. It means that they respect multiplicative structure of 
$A$, for example,
\begin{equation}
x\triangleright(ab)=(x_{(1)}\triangleright a)(x_{(2)}\triangleright b)
\,. \label{5}
\end{equation} 

One can also use this (e.g., left) covariant action of
$A^*$ upon $ A$
to construct their cross-product algebra denoted by $A\cross A^*$ 
(see, e.g.,~\cite{Ma-rev}). This is an associative 
algebra with the cross-multiplication rule given by \begin{equation} 
xa=(x_{(1)}\triangleright a)x_{(2)}\equiv \ <\!x_{(1)},a_{(2)}\!>
a_{(1)}x_{(2)}\,. \label{6}
\end{equation}
A cross-product is not a Hopf algebra but exhibits remarkable 
$A$-comodule properties~\cite{SWZ2,Sch-diss}, which
enables one to consider $A^*$ (left acting upon $A$) as an algebra 
of left-invariant (and $A\cross A^*$ of general) vector fields 
on a quantum group. 

Analogously, an associative algebra ${\cal G}=\Gamma^\wedge\cross       
\,(\Gamma^\wedge)^*$ can be defined using this cross-product
construction (here
$ (\Gamma^\wedge)^*=A^*\oplus\Gamma^*\oplus\ldots\,)$.
We place ${\cal G}$ in the center of our approach. It means that we 
assume the following guiding principle: 

{\em All cross-commutational
relations between functions, forms, Lie derivatives, and inner 
derivations are to be chosen according to the rules (\ref{6}) of a
cross-product algebra. In other words, given 
Woronowicz's calculus (and, hence, Hopf algebra 
$\Gamma^\wedge$), we then use only standard Hopf algebra technique $ 
\Gamma^\wedge \Rightarrow(\Gamma^\wedge)^*\Rightarrow \Gamma^\wedge 
\cross\,(\Gamma^\wedge)^* $ 
to construct the whole algebra of these four types of objects.} 

It only remains to put all these objects in the
corresponding `boxes'. We already know that functions and 1-forms lie
in $A$ and $\Gamma$, respectively. It is generally accepted that Lie
derivatives, ${\cal L}_h$, correspond to actions of $h\in A^*$,
\begin{equation}
{\cal L}_h\equiv h\ \triangleright\,, \label{7}
\end{equation} 
with $h\triangleright \rho $ being $\rho _{(1)}\!<\!h,\rho _{(2)}\!>$ 
for $\rho \in\Gamma^\wedge$, whereas $h\triangleright\theta 
=h_{(1)}\theta S(h_{(2)})$ \ (left adjoint action) for $\theta 
\in(\Gamma^\wedge)^*$. With this definition in mind, the covariance
(or generalized differentiation) relation (\ref{5}) is always valid.

It seems also natural to relate inner derivations with elements of
$\Gamma^*$ (cf.~\cite{Sch-Vare}). We propose the following basic 
definition:

Let $\gamma _i\in\Gamma^*$ be fixed by
\begin{equation}
<\!\gamma _i,a\omega^j\!>=\varepsilon (a)\delta _{i}^{j}\,, \label{8}
\end{equation}
and $<\!\gamma_i,\rho\!>=0$ for $\rho\in A,\Gamma^2,\Gamma^3,\ldots\,.$
Then, for $h=\chi _i$ we can define inner derivations $\iota _h$  
(denoted in this case $\iota _i$) as follows:
\begin{equation}
\iota _i\doteq\gamma _i\ \triangleright \label{9}
\end{equation}
(see~\cite{Sch-diss} for a discussion of how $\iota _h$ might be defined
for a general left-invariant field $h\in A^*$). Now, the associative 
algebra including all the four types of elements is completely 
defined.

The remaining part of the paper is devoted to drawing conclusions from 
the above definitions. Thus, one trivially shows that
\begin{equation}
\gamma _i\triangleright a=0\,, \ \ \
\gamma _i\triangleright da=\chi _i\triangleright a \,. \label{10}
\end{equation}
Not much harder (and not new, of course) is the proof that
\begin{equation}
{\cal L}_h\circ d=d\circ{\cal L}_h \,. \label{12}
\end{equation}
First, we check it for $a\in A\,: \ \ h\triangleright da=d(h
\triangleright a)\ .$ \ Then, from
$$ h\triangleright(db_1\ldots db_n)=(h_{(1)}\triangleright 
db_1)\ldots (h_{(n)}\triangleright 
db_n)=d(h_{(1)}\triangleright b_1)\ldots 
d(h_{(n)}\triangleright b_n) $$ and the Leibniz rule it follows 
that $$ h\triangleright(d(a\,db_1\ldots db_n))= 
d(h\triangleright(a\,db_1\ldots db_n))\,,  $$ which is exactly 
(\ref{12})

Applying $<\!\gamma_i,\ldots\!>$ to both $a\omega ^j$ and 
$\omega ^ja$, and then using the first eq. in
(\ref{3}), we come to
\begin{equation}
\Delta(\gamma _i)=1\otimes \gamma _i+\gamma _j\otimes f_{i}^{j}\,.
\label{13}
\end{equation} From these, 
all the multiplication relations inside ${\cal G}$,
involving $a,\omega ^i,h$ and $\gamma _i$, can be deduced more or less
automatically:
\begin{equation}
\gamma _ih=(r_{i}^{j}\triangleright h)\gamma _j\,, \label{14}
\end{equation}
\begin{equation}
<\!\gamma _i\gamma _j,a\omega ^m\omega ^n\!>=\varepsilon (a)
(\delta _{i}^{m}\delta _{j}^{n}-<\!f_{j}^{m},r_{i}^{n}\!>)\,, \label{15}
\end{equation}
and so on. However, we prefer not to list all of them here, postponing
it (as well as the detailed comparison with similar constructions of 
other authors) to a subsequent publication. Instead, we show below how
the Cartan identity
\begin{equation}
{\cal L}_h=d\circ\iota _h+\iota _h\circ d \label{16}
\end{equation}
can be derived for $h=\chi _i\,.$ One should verify that
\begin{equation}
\chi _i\triangleright\rho =d(\gamma _i\triangleright\rho )+
\gamma _i\triangleright(d\rho )\,, \ \ \ \rho \in\Gamma^\wedge\,.
\label{17}
\end{equation}

For $\rho =a\!\in\! A$ (\ref{17}) trivially follows from (\ref{10}). 
Let $\rho =a\,db \ \ (a,b\in A)\,.$ To check that
\begin{equation}
\chi _i\triangleright(a\,db)=d(\gamma _i\triangleright a\,db)+
\gamma _i\triangleright(da\wedge db)\,, \label{18}
\end{equation}
we calculate each term separately,
\begin{eqnarray*}
\chi _i\triangleright(a\,db)&=&a_{(1)}db_{(1)}<\!1\otimes \chi _i+
\chi _j\otimes f_{i}^{j}\,,\ a_{(2)}\otimes b_{(2)}\!> \\
&=&a(\chi _i\triangleright db)+(\chi _j\triangleright a)(f_{i}^{j}
\triangleright db)\,,
\end{eqnarray*}
$$ d(\gamma _i\triangleright(a\,db))=da\,(\chi _i\triangleright b)+
a\,d(\chi _i\triangleright b)\,, $$
$$ \gamma _i\triangleright(da\wedge db)=-da\,(\chi _i\triangleright b)+
(\chi _j\triangleright a)(f_{i}^{j}\triangleright db)\,, $$
and then use (\ref{12}).

At last, consider the general case $\rho =a\,db\,B\,, \ 
B=dc_1\ldots dc_n\,,
\ a,b,\ldots ,c_i\in A\,$:
$$ \chi _i\triangleright(a\,db\,B)=a\,db\,(\chi _i\triangleright B)+
a(\chi _j\triangleright db)(f_{i}^{j}\triangleright B) 
 +(\chi _k\triangleright a)(f_{j}^{k}\triangleright db)
(f_{i}^{j}\triangleright B)\,, $$
$$ \gamma _i\triangleright(a\,db\,B)=-a\,db\,(\gamma _i\triangleright 
B)+ a(\chi _j\triangleright b)(f_{i}^{j}\triangleright B)\,, $$ 
$$ 
d(\gamma _i\triangleright(a\,db\,B))=-da\,db\,(\gamma _i\triangleright 
B)+ a(\chi _j\triangleright db)(f_{i}^{j}\triangleright B) $$ $$ 
+a\,db\,d(\gamma _i\triangleright B) + da\,(\chi _j\triangleright 
b)(f_{i}^{j}\triangleright B) + a(\chi _j\triangleright 
b)\,d(f_{i}^{j}\triangleright B)\,, $$ $$ \gamma 
_i\triangleright(da\,db\,B)=da\,db\,(\gamma _i\triangleright B) 
-da\,(\chi _j\triangleright b)(f_{i}^{j}\triangleright B) $$ $$ +(\chi 
_k\triangleright a)(f_{j}^{k}\triangleright db) 
(f_{i}^{j}\triangleright B)\,. $$
After summing this up, it remains to prove that
$$ a\,db\,(\chi _i\triangleright B)=a\,db\,d(\gamma _i\triangleright 
B)\,, $$
or
$$ \chi _i\triangleright(dc_1\ldots dc_n)=d(\gamma _i\triangleright
(dc_1\ldots dc_n))\,, $$
i.e., the same problem at a lower level. Thus, the proof is completed 
by induction.

\vspace{.5cm}

{\small This work was supported in part by the INTAS grant 93-127 and 
by Russian Foundation for Basic Research (grant 95-02-05679).

I appreciate kind hospitality of M.Tarlini at the Florence section 
of INFN, where the main part of this work was done, and of 
J.Rembielinski at the University of Lodz, where it has been completed. 
I acknowledge fruitful discussions with T.Brzezinski, R.Giachetti, 
A.Isaev, P.Pyatov, J.Rembielinski, E.Sorace and M.Tarlini. I wish to 
thank O.Radko for reading the manuscript and numerous comments.}

\end{document}